\documentstyle[11pt]{article}

\title{HYPERREAL WAVES ON TRANFINITE, TERMINATED, 
DISTORTIONLESS AND LOSSLESS, TRANSMISSION LINES }
\author{A. H. Zemanian}
\date{}
\begin{document}
\maketitle
\baselineskip21pt
\newcommand{\N} {I \kern -4.5pt N}
\newcommand{\R} {I \kern -4.5pt R}

{\ Abstract --- A prior work examined the propagation of an electromagnetic 
wave on a transfinite transmission line---transfinite in the sense that 
infinitely many one-way infinite transmission lines are connected 
in cascade.  That there are infinitely many such lines results in 
the wave propagating without ever reflecting at some discontinuity
This work examines the case where the line is terminated after 
finitely many one-way infinite transmission lines with the result
that reflected waves are now produced at both the far end
as well as at the initial end of the transmisssion line.  The questions
of whether the reflected waves are infinitesimal or appreciable and whether
they sum to an infinitesimal or appreciable amount are resolved for both 
distortionless and lossless lines. Finally, the 
generalization to higher ranks of transfiniteness is briefly summarized\\

} 

\section{Introduction}

Consider a one-way-infinite electromagnetic transmission line. 
The voltage $v(x,t)$ and current $i(x,t)$ 
at the point $x$ and at time $t$ are governed by the usual wave equations
\[ \frac{\partial v(x,t)}{\partial x}\;=\;-\,\left(r\,+\,l\frac{\partial}{\partial t}\right)\,i(x,t) \]
\[ \frac{\partial i(x,t)}{\partial x}\;=\;-\,\left(g\,+\,c\frac{\partial}{\partial t}\right)\, v(x,t) \]
along with some initial and boundary conditions.  $r$, $l$, $g$, and $c$
are respectively the distributed series resistance, series inductance,
shunt conductance, and shunt capacitance per unit length of the line.
Such a conventional, one-way-infinite, transmission line will be called an 
$\omega$-{\em line}.

A natural transfinite extension of an $\omega$-line is an 
infinite cascade of $\omega$-lines, wherein the infinite extremity of each 
$\omega$-line is connected to the input of the next $\omega$-line
in the cascade, as is indicated in Figure 1.  We call such
a structure an $\omega^{2}$-{\em line}.
Of course, it takes 
an infinitely long time for any wave to pass completely through the initial 
$\omega$-line.  Furthermore, a wave will in general (but need not)
decay to an infinitesimal size when it passes beyond the initial $\omega$-line.
In any case, the wave on subsequent $\omega$-lines can 
be analyzed by means of nonstandard analysis.  In particular,
real time $t$ can be extended into unlimited hyperreal time $\bf t$, 
real distance along the initial $\omega$-line can be extended into 
unlimited hyperreal distance $\bf x$ along the $\omega^{2}$-line,
and the voltage ${\bf v}({\bf x},{\bf t})$ 
of the wave along the $\omega^{2}$-line is given as an ``internal function''
$\bf v$ of $\bf x$ and $\bf t$.  

An important property of the $\omega^{2}$-line is that it, too, is one-way 
infinite in the sense that it is a cascade of infinitely many
$\omega$-lines.  Thus, the wave front of the voltage wave
${\bf v}({\bf x},{\bf t})$ never reaches the infinite exgtremity 
of the ${\bf v}({\bf x},{\bf t})$.  No matter how large $\bf t$ is, 
that wave front will have propagated only part of the way along the 
$\omega^{2}$-line.  As a result, there is no reflected wave returning from 
that extremity.  This is the situation examined in the prior work
\cite{dist}.

On the other hand, we might think of a cascade of only finitely many 
$\omega$-lines with some terminating load resistance 
at the far end of that cascade.  We will call such a 
structure a {\em terminated transfinite line}.
It will have reflected waves.  In fact, reflections may occur at both 
the sending end $({\bf x}=0)$ and the receiving end (the far end of the line)
so that more and more waves fronts may pass forward and backward past 
any fixed point $\bf x$ of the line as $\bf t$ increases.
In general, those waves will be infinitesimally small, but not 
necessarily;  they may appreciable.  Moreover, for time $\bf t$ being a 
sufficiently large unlimited hyperreal, 
all the reflected waves will be present at a point $\bf x$, and their 
sum may or may not be infinitesimal.  Furthermore, 
these ideas also hold for terminated transmission lines having
higher ranks of transfiniteness.

The objective of this paper is to analyze such terminated transfinite lines.

This paper is written as a sequel to \cite{dist}.  We will borrow 
some results from that prior paper 
but will explain what we are borrowing.  Furthermore, we will use some
notations, terminology, and results of nonstandard analysis;
these were specified in the Introduction of \cite{dist} and will not be 
repeated here.

\section{Truncations of $\omega^{2}$-Lines and an Apparent Anomaly}

Consider the $\omega^{2}$-line of Figure 1. We wish to approximate it 
by an $\omega$-line in order to use the known results concerning waves on 
$\omega$-lines.  We do so as in \cite[Section 4]{dist} by choosing
a set of uniformly spaced sample points a distance
$\Delta x$ apart and indexed by $j=\omega k_{1}+k_{0}$.  Here, $\omega$ is 
the first transfinite ordinal.  Also, $k_{0},k_{1}\in\N$, 
where $\N$ is the set of natural numbers:
$\{0, 1, 2, \ldots\}$
For the $j$th sample point in the $\omega^{2}$-line, $k_{1}$ is the number of 
$\omega$-lines to the left of the $\omega$-line in which that $j$th 
sample point appears, and $k_{0}$ is the number of sample points to the left 
of the $j$th sample point in the $\omega$-line in 
which that $j$th sample point appears.  To obtain the $n$th 
$\omega$-line approximation of the given $\omega^{2}$-line, 
we remove the infinite part of each $\omega$-line in the $\omega^{2}$-line 
beyond the $n$th sample point of that $\omega$-line
and reconnect the resulting finite
lines in cascade.  We call this resulting
$\omega$-line the $n$th {\em truncation} of the $\omega^{2}$-line.
Thus, as $n\rightarrow \infty$, each sample point eventually 
appears in the $n$th truncation.  When this occurs, the number of 
sample points to the left of the $j$th sample point is 
$nk_{1}+k_{0}$.

With $\Delta x$ being the distance between consecutive sample points, 
the distance $x_{j,n}$ from the input to the $j$th sample point in the $n$th 
truncation is
\[x_{j,n}\;=\;(nk_{1}+k_{0})\Delta x\;\in\;\R_{+}, \]
where $\R_{+}$ denotes the nonnegative real line.

Let us henceforth assume that a voltage 
has been imposed upon the input of the 
$\omega^{2}$-line during time $t\geq 0$, with that line being 
initially at rest at $t=0$.
This induces a wave that propagates along the $\omega^{2}$-line.
With regard to the $n$th truncation, we have a wave that propagates along 
the $n$th approximating $\omega$-line, about which much analysis exists;  
see for instance \cite{we}.
When $n$ is large enough for the $j$th sample point to appear in the 
$n$th truncation, we have the voltage response at the $j$th
sample point as $v(x_{j,n},t)$.  For any positive hyperreal time 
${\bf t}=[t_{n}]>0$ and hyperreal distance ${\bf x}_{j}=
[x_{j,n}]=[(nk_{1}+k_{0})\Delta x] \geq 0$, we obtain the hyperreal voltage 
${\bf v}({\bf x}_{j},{\bf t})$ at the point ${\bf x}_{j}$
in the $\omega^{2}$-line at the time $\bf t$ specified by the 
representative sequence $\langle v(x_{j,n},t)\rangle$
indexed by $n$ for all $n$ sufficiently large.
In the usual symbolism for hyperreals, we have
\[ {\bf v}({\bf x}_{j},{\bf t})\;=\;[v(x_{j,n},t)]. \]
Note that there will be no reflected wave appearing at ${\bf x}_{j}$ 
because no 
matter how large we choose ${\bf t}=[t_{n}]$ there will be more of the
$\omega^{2}$-line further on beyond the wave front of 
${\bf v}({\bf x}_{j},{\bf t})$.  Indeed, for each $n$, the wave will have 
propagated only a finite distance along the $n$th truncation.

However, there is an apparent---but not actual---anomaly that seems
to arise now.  Consider two sample points ${\bf x}_{i}=[x_{i,n}]$
and ${\bf x}_{j}=[x_{j,n}]$, where $x_{i,n}=(nh_{1}+h_{0})\Delta x$
and $x_{j,n}=(nk_{1}+k_{0})\Delta x$.  Moreover, consider the case where 
$h_{1}<k_{1}$ but $h_{0}$ is so much larger than $k_{0}$ that 
$mh_{1}+h_{0}>mk_{1}+k_{0}$ for some positive but sufficiently small 
$m\in\N$.  Then, for $n\leq m$, $x_{j,n}$ appears in the 
$n$th truncation, while $x_{i,n}$ does not.  Consequently, 
the wave front reaches the further (from the input) sample point 
$x_{j,n}$ in the $n$th truncation (for $n\leq m$) before it appears at the 
nearer (from the input) sample point $x_{i,n}$ because the nearer one does 
not exist in that $n$th truncation.  This seems paradoxical, but it is not.
This effect occurs only for finitely many values of $n$.
For all sufficiently large $n$, both points appear
with $x_{i,n} < x_{j,n}$.  Consequently, 
the hyperreal wave front truly reaches 
${\bf x}_{i}$ before it reaches ${\bf x}_{j}$.

\section{A Terminated Transfinite Line and Its Truncations}

The ${\omega}^{2}$-line was analyzed in \cite{dist}, 
and it had no reflected waves---as 
was explained above.  Our objective in this work is
to analyze a transfinite line consisting only of finitely many 
${\omega}$-lines in cascade with perhaps a finite line as the 
last cascaded line.
As a result, a hyperreal wave along such a line will reach its 
far end at some unlimited hyperreal time $\bf T$ and 
undergo a reflection in general.  The reflected wave 
will then reach the input after another time $\bf T$ has elapsed, and another 
reflection might then occur, followed by another propagation along the line,
and so on.

More specifically, the model we shall now analyze is shown in Figure 2.
It is a transfinite line consisting of $l_{1}$ ${\omega}$-lines in cascade
with possibly a final cascaded finite line.  We index the 
1-nodes between these lines by ${\omega}$, ${\omega} 2$,
$\ldots$, ${\omega} l_{1}$.  With $k_{1}$ being 
the number of ${\omega}$-lines to the left of the ${\omega}$-line 
in which the $j$th sample point appears, we must have 
$0\leq k_{1}\leq l_{1}$.  Also, with $k_{0}$ being the number of 
sample points to the left of the $j$th sample point within the ${\omega}$-line
in which the the $j$th sample point appears, we have
$0\leq k_{0}<\infty$.  If moreover, the $j$th sample point appears in 
the last finite line, we have $0\leq k_{0}\leq l_{0}$.

The input to the line will be called the {\em sending end}, and the 
far end of the line will be called the {\em receiving end}.
Also, the {\em forward} (resp. {\em backward}) directions are from the 
sending (resp. receiving) end toward the receiving (resp. sending)
end.  Furthermore, we take it that the input is driven by a voltage source 
of value $w(t)$ in series with a sending-end resistor $R_{s}$ and a 
switch which is thrown closed at $t=0$.  So, we restrict time 
to nonnegative values.  Moreover, we assume that the line is initially
at rest at $t=0$ and that there is a resistor $R_{r}$ terminating the line at its 
receiving end.  We require that $0\leq R_{s} <\infty$ and 
$0\leq R_{r}\leq \infty$, with 0 being a possible value 
for $R_{s}$ and with 0 and $\infty$ being 
possible values for $R_{r}$.  We also assume that $0< w(t)\leq M$
for $t\geq 0$, where $M$ is a positive real number (but this 
condition on $w(t)$ can be relaxed with 
some adjustments to our conclusions).

The closing of the switch will initiate a wave that propagates forward
along the line, reflects at the receiving end, propagates backward,
reflects at the receiving end, and so on.
It will take an infinite amount of time for the wave front  
to propagate completely along any ${\omega}$-line, and the forward and backward
reflected waves might be infinitesimally small.  It is nonstandard
analysis that enables a mathematical examination of all this.

To this end, we consider a sequence of truncations of each ${\omega}$-line
into a finite line of length $n{\Delta} x$, as was done in \cite{dist}.
Let us choose ${\Delta} x$ such that the length of the last finite
line in the original terminated $\omega^{2}$-line (if such exists)
is a natural-number multiple $l_{0}$ of ${\Delta} x$.
So, for the $n$th truncation, the total length $L_{n}$ of the 
$n$th truncation of the terminated line
is $L_{n}=(nl_{1}+l_{0}){\Delta} x$.  The $j$th sample point will
appear in the $n$th truncation of the terminated $\omega^{2}$-line
when $n\geq k_{0}$ because each truncated ${\omega}$-line is of 
length $n{\Delta} x$, while $k_{0} {\Delta} x$ is the distance of the
$j$th sample point from the beginning of the ${\omega}$-line
in which that sample point appears.  Henceforth, we take it that
$n\geq k_{0}$, and thus the distance $x_{j,n}$ 
of the $j$th sample point from the sending end is 
\begin{equation}
x_{j,n}\;=\;(nk_{1}\,+\,k_{0}){\Delta} x \;=\; K_{n},  \label{3.1}
\end{equation}
and its distance $L_{n}-x_{j,n}$ to the receiving end is 
\begin{equation}
L_{n}\,-\,x_{j.n}\;=\;(n(l_{1}\,-\,k_{1})\,+\,l_{0}\,-\,k_{0}) {\Delta} x.  \label{3.2}
\end{equation}

Altogether then, each $n$th truncation is a finite transmission line
whose voltage response $v_{j}(x_{j,n},t)$ at the $j$th sample point 
at any time $t$ is uniquely determined.  So, for any hyperreal time
${\bf t}=[t_{n}]$ and with the truncated line expanding as 
$n\rightarrow\infty$ to fill out the $\omega^{2}$-line, we have 
a sequence $\langle v_{j}(x_{j,n},t_{n}\rangle$ for all $n$
sufficiently large, which determines a hyperreal response
\begin{equation}
{\bf v}({\bf x}_{j},{\bf t})\;=\;[v(x_{j,n},t_{n})]  \label{3.3}
\end{equation}
at the $j$th sample point of the $\omega^{2}$-line.  We shall now examine
that response more explicitly in terms of the reflected, 
forward and backward waves.  The latter may or 
may not be infinitesimal depending upon the parameters
of the line and the terminating resistors $R_{s}$ and $R_{r}$.
\section{A Finite Line}

Let us recall the voltage response of a finite transmission line, the
one shown in Figure 3.  In general, the propagation constant $\gamma$
is given by 
\begin{equation}
\gamma\;=\sqrt{(ls+r)(gs+c)}\;=\;\sqrt{lc}\,\sqrt{(s+\delta)^{2}-\sigma^{2}}, \;\;\;\;Re\,s\,>\,0,  \label{4.1}
\end{equation}
where $s$ is the complex variable of the Laplace transform (except 
when it appears as a subscript) and where
\begin{equation}
\delta\;=\;\frac{1}{2}\left(\frac{r}{l}+\frac{g}{c}\right),\;\;\;\;\sigma\;=\;\frac{1}{2}\left(\frac{r}{l}-\frac{g}{c}\right).  \label{4.2}
\end{equation}
Also, the characteristic impedance of the line is 
\begin{equation}
Z_{0}\;=\sqrt{\frac{ls+r}{cs+g}} .  
\end{equation}
Moreover, with $r_{s}$ and $r_{r}$ being the sending-end and receiving-end 
reflection coefficients respectively, we have
\begin{equation}
r_{s}\;=\;\frac{R_{s}-Z_{0}}{R_{s}+Z_{0}},\;\;\;\;r_{r}\;=\;\frac{R_{r}-Z_{0}}{R_{r}+Z_{0}}.  \label{4.4}
\end{equation}
We are allowing $R_{s}$ to be a short (i.e., $R_{s}=0$, $r_{s}=-1$) and
$R_{r}$ to be either a short or an open (i.e., $R_{r}=0$, $r_{r}=-1$ or
$R_{r}=\infty$, $r_{r}=1$).  With this notation and with $W(s)$ 
being the Laplace of $w(t)$, the Laplace transform $V(x,s)$
of the voltage response $v(x,t)$ at the point $x$  $(0\leq x\leq L)$
of the line at time $t\geq 0$ is the following:
\begin{eqnarray}
V(x,s)\;=\;\frac{Z_{0}}{Z_{0}+R_{s}}W(s) &(&e^{-\gamma x}\,+\,r_{r}e^{-\gamma(2L-x)} \nonumber \\
&+& r_{s}r_{r}e^{-\gamma(2L+x)}\,+\,r_{s}r_{r}^{2}e^{-\gamma(4L-x)} \nonumber \\
&+&\,\cdots \nonumber \\ 
&+&\,r_{s}^{m}r_{r}^{m}e^{-\gamma(2mL+x)}\,+\,r_{s}^{m}r_{r}^{m+1}e^{-\gamma(2(m+1)L-x)} \nonumber \\
&+&\,\cdots )\;,\;\;\;\;\;\;m=0,1,2,\ldots\; .   \label{4.5}
\end{eqnarray}
If the line's parameters $\gamma$, $R_{s}$, and $R_{r}$ are such that  
\begin{equation}
|r_{s}\,r_{r}\,e^{-2\gamma L}|\;<\; 1,  \label{4.6}
\end{equation}
this series can be summed to get 
\begin{equation}
V(x,s)\;=\;\frac{Z_{0}}{Z_{0}+R_{s}}W(s)\frac{e^{-\gamma x}\,+\,r_{r}e^{-2\gamma L}e^{\gamma x}}{1\,-\,r_{s}r_{r}e^{-2\gamma L}}.  \label{4.7}
\end{equation}

(An equation much like that of (\ref{4.5}) holds for the current on 
the line;  the analysis 
for it is similar to that for the line voltage.)

An important special case, whose transfinite version we shall examine in 
some detail, is the distortionless line.  This occurs when $r$, $g$, $l$, and $c$ are all positive
and such that 
$\sigma =0$ and $\delta=r/l=g/c$.  Thus, 
$Z_{0}$ is now the real number $\sqrt{l/c}$ and $\gamma
=\sqrt{lc}\,(s+\delta)$. 
Moreover, $-1\leq r_{s}< 1$ and $-1\leq r_{r}\leq 1$.  In this case, 
the condition (\ref{4.6}) is satisfied for Re$(s+\delta) >0$, 
and the Laplace transform (\ref{4.7})
is valid. Moreover, we may take the inverse Laplace 
transform term by term in (\ref{4.5}).  Upon setting 
$\alpha=\sqrt {lc}\,\delta$ and $u=1/\sqrt{lc}$, we obtain the inverse 
Laplace transform $v(x,t)$ of $V(x,s)$ as
\begin{eqnarray}
v(x,t)\;&=&\;\frac{Z_{0}}{Z_{0}+R_{s}}( e^{-\alpha x}w(t-x/u)\,+\,r_{r}\,e^{-\alpha(2L-x)}\,w(t-(2L-x)/u) \nonumber \\
&+&\; r_{s}\,r_{r}\,e^{-\alpha (2L+x)}\,w(t-(2L+x)/u)\,+\, r_{s}\,r_{r}^{2}\,e^{-\alpha(4L-x)}\,w(t-(4L-x)/u) \nonumber \\
&+&\; \cdots \nonumber \\
&+&\; r_{s}^{m}\,r_{r}^{m}\,e^{-\alpha(2mL+x)}\,w(t-(2mL+x)/u) \,+\, r_{s}^{m}\,r_{r}^{m+1}\,e^{-\alpha(2(m+1)L-x)}\,w(t-(2(m+1)L-x)/u) \nonumber \\
&+&\; \cdots\;).  \label{4.8} 
\end{eqnarray}
Here, $u=1/\sqrt{lc}$ is the speed of propagation of the waves.

Another special case occurs when $r=g=0$, $l>0$, and $c>0$.  Now, $\alpha=\delta =\sigma =0$.
Equation (\ref{4.8}) holds again except that all the exponential 
damping terms are replaced by 1.

\section{A Terminated Transfinite Distortionless Line}

We wish to examine the forward and backward voltage waves and also
their sum (i.e., the total voltage) at any point of the terminated
transfinite line of Figure 2, which we now take to be distortionless.
To ensure transfiniteness, we assume
that $k_{1}\geq 1$, that is, the initial line in the cascade of 
Figure 2 is truly one-way infinite.  In that initial $\omega$-line,
the wave front propagates at the speed $u=1/\sqrt{lc}$ and 
reaches a point $x$ in the $\omega$-line at the time $t=x/u$.  The voltage 
response is 0 for $t< x/u$ and possibly nonzero for $t> x/u$;  in fact;
it is positive for $t> x/u$ if $w(t)$ is positive for $t>0$.
In order to examine the propagation of the wave beyond the initial 
$\omega$-line and its subsequent reflections, we have to extend 
real time $t$ to unlimited hyperreal time ${\bf t}=[t_{n}]$, where $t_{n}
\rightarrow \infty$ as $n\rightarrow \infty$.  So, we should also extend 
the real source voltage $w$ to an internal function ${\bf w}$
of ${\bf t}$ that is identical to $w$ for all real $t$.
Specifically, we set ${\bf w}({\bf t})=[w(t_{n})]$.
We shall assume that $w(t)$ is bounded for all $t$, that is, 
there is an $M\in\R_{+}$ such that $|w(t)|\leq m$ for all $t$.
Consequently, ${\bf w}({\bf t})\leq M$ for all ${\bf t}$ as well.

To ascertain whether the forward and backward waves and their sum
are infinitesimal or instead appreciable, we examine the
$n$th truncation of the line and then determine what occurs
when $n\rightarrow \infty$.  For a more concise notation, we use 
$L_{n}=(l_{1}n+l_{0})\Delta x$ and $x_{j,n}=(k_{1}n+k_{0})\Delta x$.
Thus, $L_{n}$ is the length of the $n$th truncated line, and 
$x_{j,n}$ is the distance of the $j$th sample point from 
the sending end for all $n$ large enough to encompass
the $j$th sample point.  We always take it that $n$ is 
sufficiently large this way.  In this case, $L_{n}-x_{j,n}$ 
is the distance of the $j$th sample point from 
the receiving end.  

Thus, in accordance with (\ref{4.8}), in the $n$th truncation 
of the $\omega^{2}$-line, the voltage at $j$th sample point at time 
$t_{n}$ is obtained from (\ref{4.8}) by replacing $x$
by $x_{j,n}$, $L$ by $L_{n}$, and $t$ by $t_{n}$.
According to the first term in the resulting series,
the wave front first reaches the $j$th sample point
${\bf x}_{j}=[x_{j,n}]$ in the transfinite line
at the hyperreal time ${\bf t}={\bf x}_{j}/u=[x_{j,n}/u]$.  
Similarly, after $m$ reflections at both the receiving end 
and the sending end,
the wave front of the $m$th forward reflected wave
reaches that point at time
\[ {\bf t}\;=\;[(2mL_{n}\,+\,x_{j,n}/u]) \]
Also, after one more reflection at the receiving end,
the wave front of the $(m+1)$st backward wave reaches that point at time
\[ {\bf t}\;=\;[(2(m+1)L_{n}\,-\,x_{j,n}/u]. \]
It follows that, for ${\bf t}=[t_{n}]$ with $t_{n}=O(n)$, only finitely 
many reflected waves will have reached any hyperreal sample
point ${\bf x}_{j}$. On the other hand, if $t_{n}/n\rightarrow\infty$
as $n\rightarrow \infty$, then all of the forward and backward waves 
will have reached every hyperreal sample point ${\bf x}_{j}$.

So, now consider the hyperreal waves on the terminated transfinite 
line.
The initial forward wave in the initial $\omega$-line remains
appreciable and is given by 
\[ \frac{Z_{0}}{Z_{0}+R_{s}}[e^{-\alpha x_{j,n}}\,w(t_{n}-x_{j,n}/u)] \]
at the $j$th sample point therein, where $x_{j,n}=k_{0}\Delta x$.
However, that initial forward wave becomes infinitesimal 
beyond that initial $\omega$-line because now $x_{j,n}=
(k_{1}n+k+{0})\Delta x$ with $k_{1}>0$, 
and thus it is no larger in absolute value than 
$M\,e^{-\alpha(k_{1}n+k_{0})\Delta x}$, where $\alpha\,=\,\sqrt{lc}\,r/l\,=\, 
\sqrt{lc}\,g/c$.  Similarly, all the reflected waves are infinitesimal too.
Indeed, since $0\geq x_{j,n}\geq L_{n}$,
the absolute values of the $m$th reflected forward wave
and the $(m+1)$st reflected backward wave
are both bounded by the hyperreal 
\[ M\,[e^{-\alpha 2mL_{n}}], \] 
where $m\geq 1$, $L_{n}> 0$, and $L_{n}\rightarrow\infty$ 
as $n\rightarrow\infty$.
Thus, all the waves (i.e., the appreciable one in the initial 
$\omega$-line and all the other infinitesimal ones) diminish
at least exponentially as they propagate.

What about the sum of the reflected, forward and backward, waves?
In general, an infinite sum of infinitesimals need not
be infinitesimal, but in this case it is.  To see this, 
consider any subsequent finite line in the $n$th truncation,
the ones for which $k_{1}>0$.  (The argument is much
the same for the reflected waves in the initial line.)
Remember that $L_{n}=(l_{1}n+l_{0})\Delta x$,
$x_{j,n}=(k_{1}n+k_{0})\Delta x$, $l_{1}>0$, and
$l_{1}n+l_{0}-k_{1}n-k_{0}\geq 0$.  Also, $l_{0}\geq 0$.  So, 
\[ 2(l_{1}n+l_{0})-(k_{1}n+k_{0})\;\geq\;l_{1}n+l_{0}\;\geq l_{1}n. \] 
Since $|Z_{0}/(Z_{0}+R_{s})|\leq 1$, $|r_{s}|\leq 1$, and 
$|r_{r}|\leq 1$, we have the following estimate for the voltage 
$v(x_{j,n},t_{n})$ of the $n$th truncation at the $j$th sample point
$x_{j,n}$ beyond the initial truncated $\omega$-line:
\begin{eqnarray}
|v(x_{j,n},t_{n})|\;\leq\;M&(&e^{-\alpha(k_{1}n+k_{0})\Delta x}\;+\;e^{-\alpha(2(l_{1}n+l_{0})-(k_{1}n+k_{0}))\Delta x} \nonumber \\ 
&+&\,e^{-\alpha(2(l_{1}n+l_{0})+(k_{1}n+k_{0}))\Delta x}\;+\;e^{-\alpha(4(l_{1}n+l_{0})-(k_{1}n+k_{0}))\Delta x} \nonumber \\
&+&\,\cdots \nonumber \\
&+&\;e^{-\alpha(2m(l_{1}n+l_{0})+(k_{1}n+k_{0}))\Delta x}\;+\;e^{-\alpha(2(m+1)(l_{1}n+l_{0})-(k_{1}n+k_{0}))\Delta x} \nonumber \\
&+&\;\cdots) \nonumber
\end{eqnarray}
\begin{eqnarray}
\leq\;M&(& e^{-\alpha k_{1}n\Delta x}\;+\;e^{-\alpha l_{1}n\Delta x} \nonumber \\
       &+&\;e^{-\alpha(2l_{1}+k_{1})n\Delta x}\;+\;e^{-\alpha3l_{1}n\Delta x} \nonumber \\
       &+&\;\cdots \nonumber \\
       &+&\;e^{-\alpha(2ml_{1}+k_{1})n\Delta x}\;+\;e^{-\alpha(2m+1)l_{1}n\Delta x} \nonumber \\
       &+&\;\cdots\;\;\;\;) \nonumber
\end{eqnarray}
\[ =\;M\,\left(\frac{e^{-k_{1}n\Delta x}}{1\,-\,e^{-\alpha 2l_{1}\Delta x}}\;+\;\frac{e^{-l_{1}n\Delta x}}{1\,-\,e^{-\alpha 2l_{1}\Delta x}}\right)\,. \]
So, with regard to the $\omega^{2}$-line now,
since the right-hand side tends to 0 as $n\rightarrow\infty$, it follows 
that, at any sample point beyond the first $\omega$-line, not only each
forward or backward wave is infinitesimal but their sum is
infinitesimal, too.

This argument also shows that the reflected, forward and backward, waves
at any sample point in the initial $\omega$-line also sum to an infinitesimal.

Since there is attenuation of the wave in a distortionless line with 
$r>0$ and $g>0$, these conclusions are unsurprising at least for the 
individual, reflected, forward and backward, waves.  
Things are different for a lossless line.

\section{A Terminated Transfinite Lossless Line}

We now consider the lossless case of an $\omega^{2}$-line.  This occurs
when $r=g=0$, $l>0$, and $c>0$ 
so that $\alpha=0$, too.  Now, there is no attenuation of the
waves, but otherwise the waves propagate 
as they do in the distortionless case.
All of the forward and backward waves take on appreciable values 
at the sample points; they do not assume 
nonzero infinitesimal values.  As before, when $t_{n}=O(n)$, 
${\bf t}=[t_{n}]$, ${\bf v}({\bf x}_{j},{\bf t})$ is at most a finite sum 
of appreciable values.  If however $t_{n}/n\rightarrow\infty$ as
$n\rightarrow\infty$, ${\bf v}({\bf x}_{j},{\bf t})=[v(x_{j,n},t_{n})]$
will in general be determined by an infinite sequence of real, nonzero
values $v(x_{j,n},t_{n})$.  Just which hyperreal should be assigned to 
${\bf v}({\bf x}_{j},{\bf t})$ will depend upon 
which nonprincipal ultrafilter is chosen in our
ultrapower construction of hyperreals.

To illustrate these ideas, let us consider two examples.
As always, we assume that $k_{1}\geq 1$.

{\bf Example 6.1.}  For the sending-end source $w(t)$, we now 
choose the unit-step function $1_{+}(t)=0$ for $t<0$ and
$1_{+}(t)$ for $t\geq 0$.  Thus, each wave propagates as a unit step.
Let $0<R_{s}<\infty$ and $R_{r}=\infty$.  That is, we have a
positive resistor at the sending end and an open circuit at the 
receiving end.  So, by (\ref{4.4}), $-1<r_{s}<1$ and $r_{r}=1$.
Each forward wave is totally reflected at the receiving end
with no change of sign.  If $R_{s}\neq Z_{0}$, each backward wave 
is reflected and diminished  
by the factor $r_{s}$ (with a change of sign 
if $R_{s}< Z_{0}$) at the sending end.
There is no reflection if $R_{s}=Z_{0}$ so that $r_{s}=0$.

Altogether then, if ${\bf t}=[t_{n}]$ where $t_{n}=O(n)$, the voltage 
${\bf v}({\bf x}_{j},{\bf t})$ at the $j$th sample point 
is equal to a finite sum
\[ 2(1+r_{s}+r_{s}^{2}+r_{s}^{3}+\cdots +r_{s}^{p}). \]
On the other hand, if $t_{n}/n\rightarrow\infty$
as $n\rightarrow\infty$, ${\bf v}({\bf x}_{j},{\bf t})$ 
is infinitesimally close to 
$2/(1-r_{s})$. 

{\bf Example 6.2.}  Again, let $w(t)=1_{+}(t)$, 
but this time set $R_{s}=R_{r} =0$;  that is, $R_{s}$ and $R_{r}$ are shorts.
Thus, $r_{s}=r_{r}=-1$.  So, $r_{s}^{m}r_{r}^{m}=1$ and 
$r_{s}^{m}r_{r}^{m+1}=-1$ in (\ref{4.5}).

In fact, for any ${\bf t}=[t_{n}]$ with $t_{n}=O(n)$ and 
for any sample point
${\bf x}_{j}$, only finitely many reflected waves will have
passed ${\bf x}_{j}$.  Consequently, ${\bf v}({\bf x}_{j},{\bf t})$ 
equals either $+1$ 
or 0 depending upon whether the last reflected wave occurring at 
${\bf x}_{j}$ at time ${\bf t}$ 
is either a forward wave or a backward wave respectively.

On the other hand, if $t_{n}/n\rightarrow\infty$ as $n\rightarrow\infty$,
then the hyperreal ${\bf v}({\bf x}_{j},{\bf t})$ 
has as a representative sequence
the infinite alternating sequence: $1,0,1,0,1,0,\ldots\;$.
This means that ${\bf v}({\bf x}_{j},{\bf t})$ 
equals either 1 or 0 depending upon
the choice of the nonprincipal ultrafilter $\cal F$ used in the 
ultrapower construction of the hyperreals.  So, now we have a two-way 
ambiguity as to the value
of the hyperreal wave at any sample point ${\bf x}_{j}$ and 
time ${\bf t}$.  Without specifying $\cal F$, 
the ultrapower construction cannot 
ascertain the hyperreal value of ${\bf v}({\bf x}_{j},{\bf t})$ 
between the two possible values $1$ and 0.
\section{Transfinite, Terminated, Distortionless or Lossless, Lines of Higher Ranks}

Rather than considering a terminated $\omega^{2}$-line, as we have in the 
preceding sections, let us now consider a terminated $\omega^{\mu}$-line,
where $\mu$ is a natural number greater than 2.  The analysis of such a line
is much the same as before.  Its only essential difference is in
the more complicated expression for the terminations 
and truncations of the
$\omega^{\mu}$-line.

Let us recall the recursive definition of an $\omega^{\mu}$-line, 
where again $\mu\in\N$, $\mu>2$.  By connecting $\omega$-many 
$\omega^{2}$-lines in cascade (with the infinite extremity of each 
$\omega^{2}$-line connected to the input
of the next $\omega^{2}$-line), we obtain an $\omega^{3}$-line.
Similarly, for each $p=4,\ldots,\mu$, an $\omega^{p}$-line is
obtained by connecting $\omega$-many $\omega^{p-1}$-lines in cascade again.

As before, we choose sample points with a uniform spacing $\Delta x$
throughout the $\omega^{\mu}$-line.  Then, a typical sample point has the index
\begin{equation}
j\;=\;\omega^{\mu-1}k_{\mu-1}\,+\,\omega^{\mu-2}k_{\mu-2}\,+\,\cdots\,+\,\omega k_{1}\,+\,k_{0}  \label{7.1}
\end{equation}
where $k_{\mu-1}$ is the number of $\omega^{\mu-1}$-lines to the left of the 
$\omega^{\mu-1}$-line in which the sample point $x_{j}$ appears, $k_{\mu-2}$ 
is the number
of $\omega^{\mu-2}$-lines within the $\omega^{\mu-1}$-line in which 
the sample point 
$x_{j}$ appears and to the left of the $\omega^{\mu-2}$-line 
in which the sample 
point $x_{j}$ appears, and so on.  In general, for $1\leq\alpha\leq \mu -1$, 
$k_{\alpha}$ is the number of $\omega^{\alpha}$-lines within the 
$\omega^{\alpha+1}$-line in which the sample point $x_{j}$ appears and 
to the left of the $\omega^{\alpha}$-line in which the sample point $x_{j}$ 
appears.
Finally, $k_{0}$ is the number of sample points within the $\omega$-line
in which the sample point appears and to the left of that sample point.

Now, a terminated $\omega^{\mu}$-line is obtained by choosing 
any sample point of index
\begin{equation}
\omega^{\mu-1}l_{\mu-1}\,+\,\omega^{\mu-2}l_{\mu-2}\,+\,\cdots\,+\,\omega l_{1}\,+\,l_{0}  \label{7.1a}
\end{equation}
and deleting that part of the $\omega^{\mu}$-line beyond that point, which 
we now refer to as the receiving end of the terminated line.  Also, 
we append a receiving-end resistor $R_{r}$, as shown in Figure 2.
We now require that $l_{\mu-1}\geq 1$.  
Furthermore, we append a series connection of a voltage source $w(t)$,
a sending end resistor $R_{s}$, and a switch to be closed 
at $t=0$, again as shown in Figure 2.  
The terminated line is taken to be initially at rest, as before.

The index (\ref{7.1}) of any sample point in the terminated line
satisfies the following restriction: $0\leq k_{p}\leq l_{p}$ for 
$p=0,1,\ldots,\mu-1\;$.

The next step is to choose a truncation of the terminated 
$\omega^{\mu}$-line.  A simple way of doing this is to replace the 
``$\omega$-many'' phrase in the definition of an $\omega^{\mu}$-line
by the phrase ``$n$-many,'' where $n\in\N$.  In this case, each 
$\omega^{p}$-line $(p=1,\ldots,\mu-1)$ is replaced by an 
$n^{p}$-line and the output of each $n^{p}$-line is connected 
to the input of the next $n^{p}$-line---if there is a next
$n^{p}$-line.  Thus, the length $L_{n}$ of the terminated 
$\omega^{\mu}$-line is 
\begin{equation}
L_{n}\;=\;(n^{\mu-1}l_{\mu-1}\,+\,n^{\mu-2}l_{\mu-2}\,+\,\ldots\,+\, n l_{1}\,+\,l_{0})\,\Delta x.  \label{7.2}
\end{equation}
The distance from the sending end to any sample point in the terminated line
is
\begin{equation}
K_{n}\;=\;(n^{\mu-1}k_{\mu-1}\,+\,n^{\mu-2}k_{\mu-2}\,+\,\ldots\,+\, n k_{1}\,+\,k_{0})\,\Delta x.  \label{7.3}
\end{equation}
where $0\leq k_{p}\leq l_{p}$ $(p=0,1,\ldots, \mu-1)$, as before.
Thus, the distance from that sample point to the receiving end is 
$L_{n}-K_{n}$.

The analysis of the terminated $\omega^{\mu}$-line now proceeds 
exactly as that for an $\omega^{2}$-line except that 
$L_{n}=(l_{1}n+l_{0}) \Delta x$ is replaced by (\ref{7.2}) and 
$(k_{1}n+k_{0}) \Delta x$ is replaced by (\ref{7.3}). 
The conclusions for the distortionless line (resp. lossless line) 
are the same as those given in Sec. 5
(resp Sec. 6). 

The next stage is in this examination of terminated transfinite lines
is the examination of a terminated $\omega^{\omega}$-line.  
(See \cite[page480]{dist} for the definition of an $\omega^{\omega}$-line.)
However, such a termination yields simply a terminated
$\omega^{\mu}$-line for some natural number $\mu$, and nothing  
new transpires.

On the other hand, something more arises for an $\omega^{\mu+1}$-line.
This is a cascade of $\omega$-many $\omega^{\omega}$-lines.  
By terminating such a line at any one of its sample points, we will have 
$l_{\omega}$-many $\omega^{\omega}$-lines to the left of that sample point,
followed by an $\omega^{\mu}$-line up to the sample point.  
We can truncate this $\omega^{\mu+1}$-line, too, following the procedure 
stated in \cite[page 481]{dist} to get a finite line that expands
to fill out the $\omega^{\mu+1}$-line as the truncation index $n$
tends to $\infty$.  
Next, we can then consider an $\omega^{\omega+2}$-line, $\dots$ , an
$\omega^{\omega^{2}}$-line, and so on.  We will not pursue
these extensions, because the conclusions in every case are the same 
as those obtained heretofore.
\section{The General Case of a Transfinite, Terminated, Transmission Line}

Let us now comment very briefly about the complications that arise 
when the terminated transfinite line is neither distortionless
nor lossless.  Equation (\ref{4.5}) still holds, but now
\[ Z_{0}\;=\;\sqrt{\frac{ls+r}{cs+g}} \]
and
\[ \gamma \;=\;\sqrt{(ls+r)(cs+g)} \]
are irrational functions when $l$, $c$, and $r$ are positive and $g$ 
is nonnegative.  So, too, are the reflection coefficients $r_{s}$ 
and $r_{r}$ irrational when $0< R_{s}<\infty$ and 
$0< R_{r}<\infty$.  Thus, the terms in (\ref{4.5}) cannot 
be identified as a traveling wave as 
given by Equations (3), (4), and (5) in \cite{dist} when
$w(t)= 1_{+}(t)$. What can be concluded in this general case
remains an open question presently.

\pagebreak
\section*{Figure Captions}

\begin{description}
\item{Figure 1.} An $\omega^{2}$-line.  It consists of infinitely many 
$\omega$-lines (i.e., one-way infinite transmission lines) connected 
in cascade.  The small circles indicate
connections between the infinite extremities of the $\omega$-lines
with the inputs of the following $\omega$-lines.
\item{Figure 2.}  A terminated transfinite line; in particular,
an $(\omega l_{1}+l_{0})$-line. The number of $\omega$-lines is $l_{1}\geq 1$,
and the length of the final finite line is $l_{0} \Delta x$,
where $\Delta x$ is the distance between sample points.
\item{Figure 3.} A finite transmission line terminated at 
its sending end $(x=0)$ by a voltage source
$w(t)$ in series with a nonnegative resistor $R_{s}$ and a switch
and terminated at its receiving end $(x=L)$ by a nonnegative resistor $R_{r}$ 
(possibly $R_{r}=\infty$.)
\end{description}

\end{document}